\documentclass[12pt]{article}

\usepackage{amsmath}
\usepackage{float}
\usepackage{amsfonts}
\usepackage{amstext}
\usepackage{graphics}
\usepackage{epsf}
\usepackage[dvips]{epsfig}
\usepackage{fullpage}

\newtheorem{theorem}{Theorem}
\newtheorem{lemma}[theorem]{Lemma}
\newtheorem{corollary}[theorem]{Corollary}

\def \prend{\vrule depth-1pt height7pt width6pt}

\def \endpf{{\ \ \prend \medbreak}}

\newfont{\cyr}{wncyr10 at 12.0pt}
\newcommand{\sha}{\,\mbox{\cyr SH}\,}

\title{Avoiding Large Squares in Infinite Binary Words}

\author{Narad Rampersad, Jeffrey Shallit, and Ming-wei Wang \\
School of Computer Science \\
University of Waterloo \\
Waterloo, ON, N2L 3G1 \\
CANADA \\
{\tt nrampersad@math.uwaterloo.ca}\\
{\tt shallit@graceland.uwaterloo.ca}\\
{\tt m2wang@math.uwaterloo.ca} 
}

\begin{document}
\date{\today}
\maketitle

\begin{abstract}
We consider three aspects of avoiding large squares in infinite binary words.
First, we construct an infinite binary word avoiding both cubes $xxx$ and
squares $yy$ with $|y| \geq 4$; our construction is somewhat simpler than the original
construction of Dekking.  
Second, we construct an infinite binary word avoiding
all squares except $0^2$, $1^2$, and $(01)^2$;  our construction is somewhat simpler
than the original construction of Fraenkel and Simpson.  In both cases, we
also show how to modify our construction to obtain
exponentially many words of length $n$ with the given avoidance properties.
Finally, we answer an open
question of Prodinger and Urbanek from 1979 by demonstrating the existence of two infinite
binary words, each avoiding arbitrarily large squares, such that their perfect shuffle
has arbitrarily large squares.
\end{abstract}

\section{Introduction}

       A {\it square} is a nonempty word of the form $xx$, as in the
English word {\tt murmur}.  It is easy to see that every word of 
length $\geq 4$ constructed from the symbols $0$ and $1$ contains  
a square, so it is impossible to avoid squares in infinite binary words.
However, in 1974, Entringer, Jackson, and Schatz 
\cite{Entringer&Jackson&Schatz:1974} proved the surprising fact that
there exists an infinite binary word containing no squares $xx$ with
$|x| \geq 3$.  Further, the bound $3$ is best possible.

      A {\it cube} is a nonempty word of the
form $xxx$, as in the English sort-of-word {\tt shshsh}.
Dekking \cite{Dekking:1976} showed that there exists an infinite
binary word that contains no cubes $xxx$ and no squares
$yy$ with $|y| \geq 4$.  Furthermore, the bound $4$ is best possible.

     Dekking's construction used iterated morphisms.
By a {\it morphism} we understand a map $h: \Sigma^* \rightarrow \Delta^*$
such that $h(xy) = h(x)h(y)$ for all $x, y \in \Sigma^*$.  A morphism
may be specified by providing the {\it image words} $h(a)$ for all $a \in \Sigma$.
If $h:\Sigma^* \rightarrow \Sigma^*$ and
$h(a) = ax$ for some letter $a \in \Sigma$, then we say that
$h$ is {\it prolongable} on $a$, and we can then iterate $h$ infinitely
often to get the fixed point
$h^\omega(a) := a \, x \, h(x) \, h^2(x) \, h^3(x) \cdots $.

A morphism is $k$-{\it uniform} if $|h(a)| = k$ for all $a \in \Sigma$;
it is {\it uniform} if it is $k$-uniform for some $k$.
Uniform morphisms have particularly nice properties.  For example,
the class of words generated by iterating $k$-uniform morphisms coincides
with the class of $k$-automatic sequences, generated by finite
automata \cite{Allouche&Shallit:2003}.  

      Dekking's construction used a non-uniform morphism.  In this paper we
first show how to obtain, using the image of a uniform morphism,
an infinite binary word that is cubefree and
avoids squares $yy$ with $|y| \geq 4$.   Our construction
is somewhat simpler than Dekking's.

      Fraenkel and Simpson \cite{Fraenkel&Simpson:1995} strengthened
the results of Entringer, Jackson, and Schatz by showing that there
exists an infinite binary word avoiding all squares except
$0^2$, $1^2$, and $(01)^2$.  Their construction, however, was rather
complicated, involving several steps and non-uniform morphisms.
In this paper we show how to obtain a word where the only squares
are $0^2$, $1^2$, and $(01)^2$, using a uniform morphism.  Our construction
is somewhat simpler than that of Fraenkel and Simpson.

      We also consider the number of finite binary words satisfying the
Dekking and Fraenkel-Simpson avoidance properties.  We give exponential
upper and lower bounds on this number in both cases.

      Prodinger and Urbanek \cite{Prodinger&Urbanek:1979} also studied
words avoiding large squares, in particular with reference to operations
that preserve this property, such as the perfect shuffle $\sha$.
Let $w = a_1 a_2 \cdots a_n$ and  $x = b_1 b_2 \cdots b_n$ be words
of length $n$.  The {\it perfect shuffle} $w \sha x$ is defined to be
the word $a_1 b_1 a_2 b_2 \cdots a_n b_n$ of length $2n$.  The definition
can easily be extended to infinite words.  They stated the following open
question: do there exist two infinite words avoiding large squares such that their
perfect shuffle has arbitrarily large squares?  In this paper we resolve
this question by exhibiting an example.

\section{A cubefree word without arbitrarily long squares}

     In this section we construct an infinite cubefree binary word avoiding
squares $yy$ with $|y| \geq 4$.  The techniques we use are also used in later sections,
so in this section we spell them out in some detail.  

     We introduce the following notation for alphabets:  $\Sigma_k := \lbrace 0, 1, \ldots, k-1 \rbrace$.

\begin{theorem}
There is a squarefree infinite word over $\Sigma_4$ with
no occurrences of the subwords $12$, $13$, $21$, $32$, $231$, or $10302$.
\label{first}
\end{theorem}

\begin{proof}
Let the morphism $h$ be defined by
\begin{eqnarray*}
0 &\rightarrow& 0310201023 \\
1 &\rightarrow& 0310230102 \\
2 &\rightarrow& 0201031023 \\
3 &\rightarrow& 0203010201
\end{eqnarray*}
Then we claim the fixed point $h^\omega(0)$ has the desired properties.

First, we claim that if $w \in
\Sigma_4^*$ then $h(w)$ has no occurrences of 
$12$, $13$, $21$, $32$, $231$, or $10302$.  For if any of these
words occur as subwords of $h(w)$, they must occur within some $h(a)$ or
straddling the
boundary between $h(a)$ and $h(b)$, for some single letters $a, b$.
They do not;  this easy verification is left to the reader.

Next, we prove that if $w$ is any
squarefree word
over $\Sigma_4$ having no occurrences of 
$12$, $13$, $21$, or $32$, then $h(w)$ is squarefree.  

We argue by contradiction.  Let $w = a_1 a_2 \cdots a_n$ be a
squarefree string such that $h(w)$ contains a square,
i.e., $h(w) = xyyz$ for some $x, z \in \Sigma_4^*$,
$y \in \Sigma_4^+$.
Without loss of generality, assume that $w$ is a shortest such
string, so that $0 \leq |x|, |z| < 10$.  

Case 1:  $|y| \leq 20$.    In this case we can take $|w| \leq 5$.
To verify that $h(w)$ is squarefree,
it therefore suffices to check each of the 49 possible words $w \in
\Sigma_4^5$ to ensure that $h(w)$ is squarefree in each case.

Case 2: $|y| > 20$.  First, we establish the following result.

\begin{lemma}
\begin{itemize}
\item[(a)]
Suppose $h(ab) = t h(c) u$ for some
letters $a, b, c \in \Sigma_4$
and strings $t, u \in \Sigma_4^*$.
Then this inclusion is trivial (that is,
$t = \epsilon$ or $u = \epsilon$) or $u$ is not a prefix
of $h(d)$ for any $d \in \Sigma_4$.

\item[(b)]
Suppose there exist letters $a, b, c$ and
strings $s, t, u, v$ such that $h(a) = st$, $h(b) = uv$,
and $h(c) = sv$.  Then either $a = c$ or $b = c$.
\end{itemize}
\label{ming}
\end{lemma}

\begin{proof}
\begin{itemize}
\item[(a)]
This can be verified with a short computation.  In fact, the 
only $a, b, c$ for which the equality $h(ab) = t h(c) u$ 
holds nontrivially is $h(31) = t h(2) u$, and in this case
$t = 020301$, $u = 0102$, so $u$ is not a prefix of any $h(d)$.

\item[(b)]  This can also be verified with a short computation.
If $|s| \geq 6$, then no two distinct letters have images under $h$
that share a prefix of length $6$.
If $|s|\leq 5$, then $|t| \geq 5$, and no two distinct letters have
images under $h$ that share a suffix of length $5$.
\end{itemize}
\endpf
\end{proof}

     Once Lemma~\ref{ming} is established, the rest of the argument
is fairly standard.  It can be found, for example, in
\cite{Karhumaki&Shallit:2003}, but for completeness we repeat it here.

     For $i = 1, 2, \ldots, n$ define $A_i = h(a_i)$.  
Then if $h(w) = xyyz$, we can write
$$h(w) = A_1 A_2 \cdots A_n = A'_1 A''_1 A_2 \cdots A_{j-1} 
A'_j A''_j A_{j+1} \cdots A_{n-1} A'_n A''_n$$ 
where
\begin{eqnarray*}
A_1 &=& A'_1 A''_1 \\
A_j &=& A'_j A''_j \\
A_n &=& A'_n A''_n \\
x &=& A'_1 \\
y &=& A''_1 A_2 \cdots A_{j-1} A'_j = A''_j A_{j+1} \cdots A_{n-1} A'_n \\
z &=& A''_n, \\
\end{eqnarray*}
where $|A''_1|, |A''_j| > 0$.
See Figure~\ref{fig1}.

\begin{figure}[H]
\begin{center}
\input cube1.pstex_t
\end{center}
\caption{The string $xyyz$ within $h(w)$ \protect\label{fig1}}
\end{figure}

     If $|A''_1| > |A''_j|$, then $A_{j+1} = h(a_{j+1})$ is a subword
of $A''_1 A_2$, hence a subword of $A_1 A_2 = h(a_1 a_2)$.  Thus
we can write $A_{j+2} = A'_{j+2} A''_{j+2}$ with
$$ A''_1 A_2 = A''_j A_{j+1} A'_{j+2}.$$ 
See Figure~\ref{fig2}.

\begin{figure}[H]
\begin{center}
\input cube2.pstex_t
\end{center}
\caption{The case $|A''_1| > |A''_j|$ \protect\label{fig2}}
\end{figure}

But then,
by Lemma~\ref{ming} (a), either $|A''_j| = 0$,
or $|A''_1| = |A''_j|$, or $A'_{j+2}$ is a not
a prefix of any $h(d)$.  All three conclusions are impossible.

     If $|A''_1| < |A''_j|$, then $A_2 = h(a_2)$ is a subword of
$A''_j A_{j+1}$, hence a subword of $A_j A_{j+1} = h(a_j a_{j+1})$.
Thus we can write $A_3 = A'_3 A''_3$ with
$$ A''_1 A_2 A'_3 = A''_j A_{j+1} .$$  
See Figure~\ref{fig3}.

\begin{figure}[H]
\begin{center}
\input cube3.pstex_t
\end{center}
\caption{The case $|A''_1| < |A''_j|$ \protect\label{fig3}}
\end{figure}

By Lemma~\ref{ming} (a), either $|A''_1| = 0$ or $|A''_1| = |A''_j|$
or $A'_3$ is not a prefix of any $h(d)$.  Again, all three conclusions
are impossible.

     Therefore $|A''_1| = |A''_j|$.  
Hence $A''_1 = A''_j$, $A_2 = A_{j+1}$, $\ldots$, $A_{j-1} = A_{n-1}$,
and $A'_j = A'_n$.  Since $h$ is injective, we have
$a_2 = a_{j+1}, \ldots, a_{j-1} = a_{n-1}$.
It also follows that $|y|$ is divisible by $10$ and
$A_j = A'_j A''_j = A'_n A''_1$.   But by Lemma~\ref{ming} (b), either
(1) $a_j = a_n$ or (2) $a_j = a_1$.  In the first case,
$a_2 \cdots a_{j-1} a_j = a_{j+1} \cdots a_{n-1} a_n$, so
$w$ contains the square $(a_2 \cdots a_{j-1} a_j)^2$, a contradiction.  In the
second case, $a_1 \cdots a_{j-1} = a_j a_{j+1} \cdots a_{n-1}$, so
$w$ contains the square $(a_1 \cdots a_{j-1})^2$, a contradiction.

     It now follows that
the infinite word 
$$h^\omega(0) = 0310201023 0203010201 0310230102 0310201023 0201031023  \cdots$$
is squarefree and contains no occurrences of $12$, $13$, $21$, $32$,
$231$, or $10302$.
\endpf
\end{proof}

\begin{theorem}
Let $\bf w$ be any infinite word satisfying the conditions of Theorem~\ref{first}.
Define a morphism $g$ by
\begin{eqnarray*}
0 &\rightarrow& 010011 \\
1 &\rightarrow& 010110 \\
2 &\rightarrow& 011001 \\
3 &\rightarrow& 011010
\end{eqnarray*}
Then $g({\bf w})$ is a cubefree word containing no squares $xx$ with
$|x| \geq 4$.
\label{second}
\end{theorem}

     Before we begin the proof, we remark that all the words
$12$, $13$, $21$, $32$, $231$, $10302$ must indeed be avoided,
because
\begin{eqnarray*}
g(12) &&\mbox{contains the squares}\ (0110)^2,\ (1100)^2,\ (1001)^2 \\
g(13) &&\mbox{contains the square}\ (0110)^2 \\
g(21) &&\mbox{contains the cube}\ (01)^3 \\
g(32) &&\mbox{contains the square}\ (1001)^2 \\
g(231)  &&\mbox{contains the square}\ (10010110)^2 \\
g(10302) &&\mbox{contains the square}\ (100100110110)^2.
\end{eqnarray*}

\begin{proof}
     The proof parallels the proof of Theorem~\ref{first}.   
     Let $w = a_1
a_2 \cdots a_n$ be a squarefree string, with no occurrences of 
$12$, $13$, $21$, $32$, $231$, or $10302$.  We first establish that if
$g(w) = xyyz$ for some
$x, z \in \Sigma_4^*$, $y \in \Sigma_4^+$,
then $|y| \leq 3$.
Without loss of generality, assume $w$ is a shortest such string, so
$0 \leq |x|, |z| < 6$.

     Case 1:  $|y| \leq 12$.  In this case we can take $|w| \leq 5$.  To verify
that $g(w)$ contains no squares $yy$ with
$|y| \geq 4$, it suffices
to check each of the $41$ possible words $w \in \Sigma_4^5$.

      Case 2:  $|y| > 12$.  First, we establish the analogue of
Lemma~\ref{ming}.

\begin{lemma}
\begin{itemize}
\item[(a)]
Suppose $g(ab) = t g(c) u$ for some
letters $a, b, c \in \Sigma_4$
and strings $t, u \in \Sigma_4^*$.
Then this inclusion is trivial (that is,
$t = \epsilon$ or $u = \epsilon$) or $u$ is not a prefix
of $g(d)$ for any $d \in \Sigma_4$.

\item[(b)]
Suppose there exist letters $a, b, c$ and
strings $s, t, u, v$ such that $g(a) = st$, $g(b) = uv$,
and $g(c) = sv$.  Then either $a = c$ or $b = c$, or
$a = 2$, $b = 1$, $c = 3$, $s = 0110$, $t = 01$, $u = 0101$, $v = 10$.
\end{itemize}
\label{ming2}
\end{lemma}

\begin{proof}
\begin{itemize}
\item[(a)]  This can be verified with a short computation.  
The only $a, b, c$ for which $g(ab) = t g(c) u$ holds nontrivially
are 
\begin{eqnarray*}
g(01) &=& 010 \ g(3) \ 110 \\
g(10) &=& 01 \ g(2) \ 0011 \\
g(23) &=& 0110 \ g(1) \ 10 .\\
\end{eqnarray*}
But none of $110$, $0011$, $10$ are prefixes of any $g(d)$.

\item[(b)] If $|s| \geq 5$ then no two distinct letters have images under
$g$ that share a prefix
of length $5$.  If $|s| \leq 3$ then $|t| \geq 3$, and no two distinct
letters share a suffix of length $3$.  Hence $|s| = 4$, $|t| = 2$.
But only $g(2)$ and $g(3)$ share a prefix of length $4$, and only
$g(1)$ and $g(3)$ share a suffix of length $2$.
\end{itemize}
\endpf
\end{proof}

     The rest of the proof is exactly parallel to the proof
of Theorem~\ref{first}, with the following exception.  When we
get to the final case, where $|y|$ is divisible by $6$, we can use
Lemma~\ref{ming2} to
rule out every case except where $x = 0101$, $z = 01$,
$a_1 = 1$, $a_j = 3$, and $a_n = 2$.  Thus $w = 1 \alpha 3 \alpha 2$
for some string $\alpha \in \Sigma_4^*$.  This
special case is ruled out by the following lemma:

\begin{lemma}
     Suppose $\alpha \in \Sigma_4^*$, and let
$w = 1 \alpha 3 \alpha 2$.  Then either $w$ contains a square, or
$w$ contains an occurrence of one of the subwords $12$, $13$,
$21$, $32$, $231$, or $10302$.
\end{lemma}

\begin{proof}
      This can be verified by checking (a) all strings $w$ 
with $|w| \leq 4$, and (b) all strings of the form $w = abc w' de$, where
$a, b, c, d, e \in \Sigma_4$ and $w' \in
\Sigma_4^*$.  (Here $w'$ may be treated as an
indeterminate.)
\endpf
\end{proof}

    It now remains to show that if $w$ is squarefree and contains
no occurrence of $12$, $13$, $21$, $32$, $231$, or $10302$, then
$g(w)$ is cubefree.  If $g(w)$ contains 
a cube $yyy$, then it contains a square $yy$, and from what precedes
we know $|y| \leq 3$.  It therefore suffices to show that $g(w)$
contains no occurrence of $0^3$, $1^3$, $(01)^3$, $(10)^3$,
$(001)^3$, $(010)^3$, $(011)^3$, $(100)^3$, $(101)^3$, $(110)^3$.
The longest such string is of length $9$, so it suffices to examine the
$16$ possibilities for $g(w)$ where $|w| = 3$.  This is left to the reader.

    The proof of Theorem~\ref{second} is now complete.
\endpf
\end{proof}

\begin{corollary}
    If $g$ and $h$ are defined as above, then
$$g(h^\omega(0)) = 
010011011010010110010011011001010011010110010011011001011010
\cdots$$
is cubefree, and avoids all squares $xx$ with $|x| \geq 4$.
\label{jcor}
\end{corollary}

    Next, based on the morphism $h$, we define the substitution
$h':\Sigma_4^* \rightarrow 2^{\Sigma_4^*}$ as follows:
\begin{eqnarray*}
0 &\rightarrow& \{h(0)\} \\
1 &\rightarrow& \{h(1), 0310230201\} \\
2 &\rightarrow& \{h(2)\} \\
3 &\rightarrow& \{h(3)\}
\end{eqnarray*}
Thus, if $w\in\Sigma_4^*$, $h'(w)$ is a language of $2^r$ words over
$\Sigma_4$, where $r=|w|_1$.  Each of these words is of length $10|w|$.

\begin{lemma}
    Let $g$, $h$, and $h'$ be defined as above.  Let $w=h^m(0)$ for some
positive integer $m$.  Then $g(h'(w))$ is a language of $2^{n/300}$
words over $\Sigma_2$, where $n=60\cdot 10^m$ is the length of each
of these words.  Furthermore, each of these words is cubefree and avoids all 
squares $xx$ with $|x|\geq 4$.
\label{lowbnd1}
\end{lemma}

\begin{proof}
    Note that there are exactly two 1's in every image word of $h$.
Hence, $|w|_1 = {1 \over 5}|w| = {1 \over 5}\cdot 10^m$.  We have then that
$g(h'(w))$ consists of $2^{{1 \over 5}\cdot 10^m}$ binary words.  Since
$n=6\cdot 10\cdot 10^m$, we see that $g(h'(w))$ consists of $2^{n/300}$
words.

    To see that the words in $g(h'(w))$ are cubefree and avoid all
squares $xx$ with $|x|\geq 4$, it suffices by Theorem~\ref{second}
to show that the words in $h'(w)$ are squarefree and contain
no occurrences of the subwords 12, 13, 21, 32, 231, or 10302.
By the same reasoning as in Theorem~\ref{first}, the reader may easily
verify that no word in $h'(w)$ contains an occurrence of 12, 13, 21, 32,
231, or 10302.

    To show that the words in $h'(w)$ are squarefree, we will, as a
notational convenience, prefer to consider $h'$ to be a morphism defined
as follows:
\begin{eqnarray*}
0 &\rightarrow& h(0) \\
1 &\rightarrow& h(1) \\
\hat{1} &\rightarrow& 0310230201 \\
2 &\rightarrow& h(2) \\
3 &\rightarrow& h(3)
\end{eqnarray*}
Here, 1 and $\hat{1}$ are considered to be the same alphabet symbol; the `hat'
simply serves to distinguish between which choice is made for the
substitution.  To show that $h'(w)$ is squarefree, it suffices to
show that $h'$ satisfies the conditions of Lemma~\ref{ming}.  For
Lemma~\ref{ming}~(a) we have $h'(22) = th'(\hat{1})u$, but we can rule this
case out since $w$ avoids the square 22.  For Lemma~\ref{ming}~(b)
we again have that no two distinct letters have images under $h'$ that
share a prefix of length 6 or a suffix of length 5
(since 1 and $\hat{1}$ are not considered to be distinct letters).
Hence, $h'$ satisfies the conditions of Lemma~\ref{ming}, and so
$h'(w)$ is squarefree.
\endpf
\end{proof}

\begin{theorem}
    Let $G_n$ denote the number of cubefree binary words of length $n$ that
avoid all squares $xx$ with $|x|\geq 4$.  Then $G_n = \Omega(1.002^n)$ and
$G_n = O(1.178^n)$.
\label{numwords1}
\end{theorem}

\begin{proof}
    Noting that $2^{1/300} \doteq 1.002$, we see that the lower bound follows
immediately from Lemma~\ref{lowbnd1}.

    For the upper bound we reason as follows.  The set of binary words of
length $n$ avoiding cubes and squares $xx$ with $|x|\geq 4$ is a a subset
of the set of binary words avoiding 000 and 111.  The number $G'_n$
of binary words avoiding 000 and 111 satisfies the linear recurrence
$G'_n = G'_{n-1} + G'_{n-2}$ for $n\geq 3$.  From well-known properties
of linear recurrences, it follows that $G'_n = O(\alpha^n)$, where
$\alpha$ is the largest zero of $x^2 - x - 1$, the characteristic
polynomial of the recurrence.  Here $\alpha < 1.618$, so $G'_n = O(1.618^n)$.

    This argument can be extended by using a symbolic algebra package such
as Maple.  Noonan and Zeilberger \cite{Noonan&Zeilberger:1999}
have written a Maple package, \verb=DAVID_IAN=, that allows one to specify
a list $L$ of forbidden words, and computes the generating function
enumerating words avoiding members of $L$.  We used the package
for a list $L$ of 90 words of length $\leq 20$:
$$000,111,\ldots,11011001001101100100$$
obtaining a characteristic polynomial of degree 44 with dominant root
$\doteq 1.178$.
\endpf
\end{proof}

     The following table gives the number $G_n$ of binary words of length
$n$ avoiding both cubes $xxx$ and squares $y$ with $|y| \geq 4$.
\begin{center}
\begin{tabular}{|c|rrrrrrrrrrrrrrrrrr|}
\hline
$n$& 0 & 1 & 2 & 3 & 4 & 5 & 6 & 7 & 8 & 9 & 10 & 11 & 12 & 13 & 14 & 15 & 16 & 17 \\
\hline
$G_n$& 1 & 2 & 4 & 6 & 10 & 16 & 24 & 36 & 52 & 72 & 90 & 116 & 142 & 178 & 220 & 264
& 332 & 414 \\
\hline
\end{tabular}
\end{center}

\section{A uniform version of Fraenkel-Simpson}

     In this section we construct an infinite binary word avoiding all squares
except $0^2$, $1^2$, and $(01)^2$.

       Roughly speaking, verifying that the image of a morphism
avoids arbitrarily large squares breaks up into two
parts:  checking a finite number of ``small'' squares, and checking an
infinite number of ``large'' squares.  The small squares can be checked
by brute force, while for the large squares we need a version of
Lemma~\ref{ming}.  Referring to Lemma~\ref{ming} (a),
if $h(c)$ is a subword of $h(ab)$ for some letters
$a, b, c$, we call this an ``inclusion''.  Inclusions can be ruled out either
by considering prefixes, as we did in Lemma~\ref{ming} (a), or suffixes.
Referring to Lemma~\ref{ming} (b),
if $h(a) = st$, $h(b) = uv$, and $h(c) = sv$, we call that an ``interchange''.

The basic idea of the proofs in this section parallels that of the previous
section, so we just sketch the basic ideas, pointing out the properties of the
inclusions and interchanges.

     Consider the 24-uniform morphism $h$ defined as follows:
\begin{eqnarray*}
0 &\rightarrow&  012321012340121012321234 \\
1 &\rightarrow&  012101234323401234321234 \\
2 &\rightarrow&  012101232123401232101234 \\
3 &\rightarrow&  012321234323401232101234 \\
4 &\rightarrow&  012321234012101234321234
\end{eqnarray*}

\begin{theorem}
If $w \in \Sigma_5^*$ is squarefree and avoids the patterns
$02, 03, 04, 14, 20, 30, 41$, then $h(w)$ is squarefree and
avoids the patterns
$02, 03, 04, 13, 14, 20, 24, 30, 31, 41, 42, 010, 434$.
\label{third}
\end{theorem}

\begin{proof}
       The only inclusion is
$h(32) = 0123212343234 \ h(0) \  01232101234$, and
$0123212343234$ is not a suffix of the image of any letter.

     There are no interchanges for this morphism.
\endpf
\end{proof}

Now consider the 6-uniform morphism

\begin{eqnarray*}
g(0) &=& 011100 \\
g(1) &=& 101100 \\
g(2) &=& 111000 \\
g(3) &=& 110010 \\
g(4) &=& 110001
\end{eqnarray*}

\begin{theorem}
If $w$ is squarefree and avoids the patterns
$02, 03, 04, 13, 14, 20, 24, 30, 31, 41, 42, 434010$
then the only squares in $g(w)$ are $00, 11, 0101$.
\label{fourth}
\end{theorem}

\begin{proof}

There are no examples of interchanges for $g$.

There are multiple examples of inclusions, but many of them can be ruled
out by properties of $w$ and $g$:

\begin{itemize}
\item $g(02) = 01110 g(0) 0$  but $02$ cannot occur
\item $g(24) = 1 g(4) 10001$  but $24$ cannot occur
\item $g(12) = 10110 g(0) 0$  but $10110$ is not a suffix of any $g(a)$
\item $g(32) = 11001 g(0) 0$  but $11001$ is not a suffix of any $g(a)$
\item $g(21) = 1 g(4) 01100$  but $01100$ is not a prefix of any $g(a)$
\item $g(23) = 1 g(4) 10010$  but $10010$ is not a prefix of any $g(a)$
\end{itemize}

Since $g(434010) = 1100 (01110010110001)^2 1100$,
we need a special argument to rule this out.
There are four special cases that must be handled:

\begin{itemize}
\item $g(43) = 1100 g(0) 10$
\item $g(34) = 1100 g(1) 01$
\item $g(01) = 01 g(3) 1100$
\item $g(10) = 10 g(4) 1100$
\end{itemize}

In the first example, $g(43) = 1100 g(0) 10$, since $10$ is only a prefix of $g(1)$,
we can extend on the right to get $g(43) 1100 = 1100 g(01)$.  But since
$1100$ is only a prefix of $g(3)$ or $g(4)$, this gives either the forbidden
pattern $33$ or the forbidden pattern $434$.

In the second example, $g(34) = 1100 g(1) 01$, since $01$ is only a prefix of $g(0)$,
we can extend on the right to get $g(34) 1100 = 1100 g(10)$.  
But $1100$ is a suffix of only $g(0)$ and $g(1)$, so on the right we get either
the forbidden pattern $010$ or the forbidden pattern $11$.

The other two cases are handled similarly.  \endpf

%
%
%
%
\end{proof}

    As in the previous section, we now define the substitution
$h':\Sigma_5^* \rightarrow 2^{\Sigma_5^*}$ as follows:
\begin{eqnarray*}
0 &\rightarrow&  \{h(0), 012101232123401234321234\} \\
1 &\rightarrow&  \{h(1)\} \\
2 &\rightarrow&  \{h(2)\} \\
3 &\rightarrow&  \{h(3)\} \\
4 &\rightarrow&  \{h(4)\}
\end{eqnarray*}
Thus, if $w\in\Sigma_5^*$, $h'(w)$ is a language of $2^r$ words over
$\Sigma_5$, where $r=|w|_0$.  Each of these words is of length $24|w|$.

\begin{lemma}
    Let $g$, $h$, and $h'$ be defined as above.  Let $w=h^m(0)$ for some
positive integer $m$.  Then $g(h'(w))$ is a language of $2^{n/1152}$
words over $\Sigma_2$, where $n=144 \cdot 24^m$ is the length of each
of these words.  Furthermore, these words avoid all squares except
$0^2$, $1^2$, and $(01)^2$.
\label{lowbnd2}
\end{lemma}

\begin{proof}
    The proof is analogous to that of Lemma~\ref{lowbnd1}.
Note that there are at least three 0's in every image word of $h$.
Hence, $|w|_0 \geq {1 \over 8}|w| = {1 \over 8}\cdot 24^m$.  We have then that
$g(h'(w))$ consists of at least $2^{{1 \over 8}\cdot 24^m}$ binary words.  Since
$n=6\cdot 24\cdot 24^m$, we see that $g(h'(w))$ consists of at least
$2^{n/1152}$ words.

    To see that the words in $g(h'(w))$ avoid all squares except
$0^2$, $1^2$, and $(01)^2$ it suffices by Theorem~\ref{fourth}
to show that the words in $h'(w)$ are squarefree and contain
no occurrences of the subwords 02, 03, 04, 13, 14, 20, 24, 30, 31, 41, 42,
or 434010.  The reader may easily verify that the words in $h'(w)$
contain no occurrences of the subwords 02, 03, 04, 13, 14, 20, 24, 30, 31, 41,
42, or 434010.

    To show that the words in $h'(w)$ are squarefree, we will, as
before, consider $h'$ to be a morphism defined as follows:
\begin{eqnarray*}
0 &\rightarrow& h(0) \\
\hat{0} &\rightarrow& 012101232123401234321234 \\
1 &\rightarrow& h(1) \\
2 &\rightarrow& h(2) \\
3 &\rightarrow& h(3) \\
4 &\rightarrow& h(4)
\end{eqnarray*}
There are no inclusions for $h'$ other than the one identified in the
proof of Theorem~\ref{third}.  There are three interchanges: referring
to Lemma~\ref{ming}~(b), we have that $(a,b,c) \in
\{(2,1,\hat{0}),(2,4,\hat{0}), (\hat{0},3,2)\}$ satisfies $h'(a) = st$, $h'(b) = uv$,
and $h'(c) = sv$.  We may rule out the first two cases by showing that $w$ avoids
all subwords of the form $1\alpha0\alpha2$ and $4\alpha0\alpha2$,
where $\alpha \in \Sigma_5^*$.
Note that in the word $w$, any occurrence of 0 must be followed by a 1,
since $w$ avoids the patterns 02, 03, and 04.  Let $x$ be a subword of $w$
of the form $1\alpha0\alpha2$ or $4\alpha0\alpha2$.  Then $x$ must begin
with 11 or 41.  This is a contradiction, as $w$ avoids both 11 and 41.

We may rule out the third case by showing that $w$
avoids all subwords of the form $3\alpha2\alpha0$, where
$\alpha\in\Sigma_5^*$.  Note that in the word $w$, any occurrence of 2
must be followed by either 1 or 3, since $w$ avoids the patterns 20 and
24.  Let $x$ be a subword of $w$ of the form $3\alpha2\alpha0$.  Then $x$
must begin with 31 or 33.  This is a contradiction, as $w$ avoids both 31
and 33.
\endpf
\end{proof}

\begin{theorem}
    Let $H_n$ denote the number of binary words of length $n$ that
avoid all squares except $0^2$, $1^2$, and $(01)^2$.  Then
$H_n = \Omega(1.0006^n)$ and $H_n = O(1.135^n)$.
\label{numwords2}
\end{theorem}

\begin{proof}
    The proof is analogous to that of Theorem~\ref{numwords1}.
Noting that $2^{1/1152} \doteq 1.0006$, we see that the lower bound follows
immediately from Lemma~\ref{lowbnd2}.

    For the upper bound, we again used the \verb=DAVID_IAN= Maple package
for a list of 65 words of length $\leq 20$:
$$0000,1010,\ldots,1110001011100010$$
obtaining a characteristic polynomial of degree 58 with dominant root
$\doteq 1.135$.
\endpf
\end{proof}

     The following table gives the number $H_n$ of binary words of length
$n$ containing only the squares $0^2$, $1^2$, and $(01)^2$.
\begin{center}
\begin{tabular}{|c|rrrrrrrrrrrrrrrrrrr|}
\hline
$n$& 0 & 1 & 2 & 3 & 4 & 5 & 6 & 7 & 8 & 9 & 10 & 11 & 12 & 13 & 14 & 15 & 16 & 17 & 18\\
\hline
$H_n$& 1 & 2 & 4 & 8 & 13 & 22 & 31 & 46 & 58 & 78 & 99 & 124 & 144 & 176 & 198 & 234
& 262 & 300 & 351 \\
\hline
\end{tabular}
\end{center}

\section{The Prodinger-Urbanek problem}
\label{pu}

       Prodinger and Urbanek \cite{Prodinger&Urbanek:1979} stated they were
unable to find an example of two infinite binary words
avoiding large squares such that their perfect shuffle had arbitrarily large squares.
In this section we give an example of such words.

\begin{theorem}
There exist two infinite binary words $\bf x$ and $\bf y$ such that neither $\bf x$
nor $\bf y$ contain a square $ww$ with $|w| \geq 4$, but ${\bf x} \sha {\bf y}$
contains arbitrarily large squares.
\label{pu-thm}
\end{theorem}

\begin{proof}
Consider the morphism $f:\Sigma_2^*\rightarrow\Sigma_2^*$ defined as follows:
\begin{eqnarray*}
f(0) & = & 001 \\
f(1) & = & 110.
\end{eqnarray*}
We will show that
$f^\omega(0) = 001001110001001110110110001 \cdots $ contains arbitrarily large
squares and is the perfect shuffle of two words,
each avoiding squares $ww$ with $|w| \geq 4$.

First, we define the morphisms $h:\Sigma_4^*\rightarrow\Sigma_4^*$,
$g_1:\Sigma_4^*\rightarrow\Sigma_2^*$, and
$g_2:\Sigma_4^*\rightarrow\Sigma_2^*$ defined as follows:
\begin{eqnarray*}
h(0) & = & 012 \\
h(1) & = & 302 \\
h(2) & = & 031 \\
h(3) & = & 321,
\end{eqnarray*}
\vspace{-5ex}

\begin{minipage}{2.5in}
\begin{eqnarray*}
g_1(0) & = & 001 \\
g_1(1) & = & 101 \\
g_1(2) & = & 010 \\
g_1(3) & = & 110,
\end{eqnarray*}
\end{minipage}
\hfill\parbox{0.5in}{and}\hfill
\begin{minipage}{2.5in}
\begin{eqnarray*}
g_2(0) & = & 010 \\
g_2(1) & = & 100 \\
g_2(2) & = & 011 \\
g_2(3) & = & 101.
\end{eqnarray*}
\end{minipage}
\bigskip

We now show
\begin{lemma}
\label{Shuf}
$f^\omega(0) = g_2(h^\omega(0)) \sha g_1(h^\omega(0)).$
\end{lemma}

\begin{proof}
We prove the following identities by induction on $n$.
\begin{eqnarray}
\label{00}
f^{n+1}(00) & = & g_2(h^n(0)) \sha g_1(h^n(0)) \\
f^{n+1}(10) & = & g_2(h^n(1)) \sha g_1(h^n(1)) \\
f^{n+1}(01) & = & g_2(h^n(2)) \sha g_1(h^n(2)) \\
f^{n+1}(11) & = & g_2(h^n(3)) \sha g_1(h^n(3))
\end{eqnarray}
It is easy to verify that these equations hold for $n=0$.
We assume that they hold for $n=k$, where $k>0$, and show that they hold
for $n=k+1$.  We first consider $f^{k+2}(00)$, where we have
\begin{eqnarray*}
f^{k+2}(00) & = & f^{k+1}(001001) \\
& = & f^{k+1}(00)f^{k+1}(10)f^{k+1}(01) \\
& = & \left(g_2(h^k(0))\sha\,g_1(h^k(0))\right)\,
      \left(g_2(h^k(1))\sha\,g_1(h^k(1))\right)\,
      \left(g_2(h^k(2))\sha\,g_1(h^k(2))\right) \\
& = & \left(g_2(h^k(0))g_2(h^k(1))g_2(h^k(2))\right)\sha\,
      \left(g_1(h^k(0))g_1(h^k(1))g_1(h^k(2))\right) \\
& = & g_2(h^k(0)h^k(1)h^k(2))\sha\,g_1(h^k(0)h^k(1)h^k(2)) \\
& = & g_2(h^k(012))\sha\,g_1(h^k(012)) \\
& = & g_2(h^{k+1}(0))\sha\,g_1(h^{k+1}(0))
\end{eqnarray*}
as desired.  The other cases of the induction for $f^{k+2}(10), f^{k+2}(01)$,
and $f^{k+2}(11)$ follow similarly.  The result now follows from \eqref{00}.
\endpf
\end{proof}

We now prove 
\begin{lemma}
The infinite word $h^\omega(0)$ is squarefree.
\label{hSQF}
\end{lemma}

\begin{proof}
This follows immediately by the analogue of Lemma~\ref{ming}.  An easy
computation shows there are no inclusions or interchanges for $h$.
\endpf
\end{proof}

We now define
$$\mathcal{A}=\{010, 013, 021,
030, 032, 102, 121, 131, 202, 212, 231, 301, 303, 312, 320, 323\}.$$

\begin{lemma}
\label{hAvoid}
\begin{enumerate}
\renewcommand{\labelenumi}{(\alph{enumi})}
\item $h^\omega(0)$ contains no subwords $x$ where $x\in \mathcal{A}$; and

\item $h^\omega(0)$ contains no subwords of the form $0\alpha1\alpha3$,
$1\alpha0\alpha2$, $2\alpha3\alpha1$, or $3\alpha2\alpha0$, where
$\alpha\in\Sigma_4^*$.
\end{enumerate}
\end{lemma}

\begin{proof}
\begin{enumerate}
\renewcommand{\labelenumi}{(\alph{enumi})}
\item This can be verified by inspection.

\item We argue by contradiction.  Let $w$ be a shortest subword of
$h^\omega(0)$ such that $w$ is of the form $0\alpha1\alpha3$,
$1\alpha0\alpha2$, $2\alpha3\alpha1$, or $3\alpha2\alpha0$.
Suppose $w$ is of the form $0\alpha1\alpha3$. Note that the only
image words of $h$ that contain the letter 1 are $h(0)=012$,
$h(2)=031$, and $h(3)=321$.  Hence it must be the case that $\alpha$
is of the form $2\alpha'0$, $\alpha'03$, or $\alpha'32$ for some
$\alpha'\in\Sigma_4^*$.  We therefore have three cases.

\begin{enumerate}
\renewcommand{\labelenumii}{Case \arabic{enumii}:}
\item $w=02\alpha'012\alpha'03$ for some $\alpha'\in\Sigma_4^*$.
We have two subcases.

\begin{enumerate}
\renewcommand{\labelenumiii}{Case \arabic{enumii}.\roman{enumiii}:}
\item $|w|\leq 12$.  A short computation suffices to verify
that, contrary to (a), all words $w$ of the form $02\alpha'012\alpha'03$
with $|w|\leq 12$ contain a subword $x$ where $x\in\mathcal{A}$.

\item $|w|>12$.  We first make the observation that any image word of
$h$ is uniquely specified by its first two letters and also by its
last two letters.  Thus, if $w=02\alpha'012\alpha'03$, it must be the
case that $3w1=302\alpha'012\alpha'031=h(1)\alpha'h(0)\alpha'h(2)$
is also a subword of $h^\omega(0)$.  Furthermore, since $h$ has no inclusions, the
infinite word $h^\omega(0)$ can be uniquely parsed into image words
of $h$.  Since we have that $h(1)\alpha'h(0)\alpha'h(2)$ is a subword of
$h^\omega(0)$, this implies that $|\alpha'|$ is a multiple of 3 and
that $h(1)\alpha'h(0)\alpha'h(2)=h(1\beta0\beta2)$
for some $\beta\in\Sigma_4^*$, $|\beta|<|\alpha'|$.
So $1\beta0\beta2$ must also be a subword of $h^\omega(0)$.
This contradicts the minimality of $w$.
\end{enumerate}

\item $w=0\alpha'031\alpha'033$ for some $\alpha'\in\Sigma_4^*$.  But
then $w$ contains the square 33, contrary to Lemma~\ref{hSQF}.

\item $w=0\alpha'321\alpha'323$ for some $\alpha'\in\Sigma_4^*$.  But
by (a) $w$ cannot contain the subword 323.
\end{enumerate}

The cases where $w$ is of the form 
$1\alpha0\alpha2$, $2\alpha3\alpha1$, or $3\alpha0\alpha2$
follow similarly.
\end{enumerate}
\endpf
\end{proof}

We now give the analogue of Lemma~\ref{ming} for $g_1$ and $g_2$.
Let $g_i$ represent either $g_1$ or $g_2$.  Then we have
\begin{lemma}
\label{gMainLem}
\begin{enumerate}
\renewcommand{\labelenumi}{(\alph{enumi})}
\item Suppose $g_i(ab)=tg_i(c)u$ for some letters $a,b,c\in\Sigma_4$ and
words $t,u\in\Sigma_2^*$.  Then at least one of the following holds:
\begin{enumerate}
\renewcommand{\labelenumii}{(\roman{enumii})}
\item this inclusion is trivial (that is, $t=\epsilon$ or $u=\epsilon$);
\item $u$ is not a prefix of $g_i(d)$ for any $d\in\Sigma_4$;
\item $t$ is not a suffix of $g_i(d)$ for any $d\in\Sigma_4$; or
\item for all $v,w\in\Sigma_2^*$ and all $e,e'\in\Sigma_4$, if
$vg_i(ab)w=g_i(ece')$, then at least one of the following holds:
\begin{enumerate}
\renewcommand{\labelenumiii}{(\Alph{enumiii})}
\item this inclusion is trivial (that is, $v=\epsilon$ or $w=\epsilon$);
\item $w$ is not a prefix of $g_i(d)$ for any $d\in\Sigma_4$;
\item $v$ is not a suffix of $g_i(d)$ for any $d\in\Sigma_4$;
\item either $e=c$ or $e'=c$;
\item $ece'\in\mathcal{A}$;
\item for all $x,y\in\Sigma_2$ and all $k\in\Sigma_4$, if
$g_i(kab)x=yg_i(ece')$, then $k=a$; or
\item for all $x,y\in\Sigma_2$ and all $k\in\Sigma_4$, if
$xg_i(abk)=g_i(ece')y$, then $k=b$.
\end{enumerate}
\end{enumerate}

\item Suppose there exist letters $a,b,c\in\Sigma_4$ and words
$s,t,u,v\in\Sigma_2^*$ such that $g_i(a)=st$, $g_i(b)=uv$, $g_i(c)=sv$,
and $b\alpha c\alpha a$ is a subword of $h^\omega(0)$ for some
$\alpha\in\Sigma_4^*$.  Then either $a=c$ or $b=c$.
\end{enumerate}
\end{lemma}

\begin{proof}
\begin{enumerate}
\renewcommand{\labelenumi}{(\alph{enumi})}
\item We give one example of each case and list the other non-trivial
cases in a table below.  
\begin{enumerate}
\renewcommand{\labelenumii}{(\roman{enumii})}
\item Trivial.
\item $g_2(32)=1g_2(0)11$, but 11 is not a prefix of $g_2(d)$ for any
$d\in\Sigma_4$.
\item $g_1(02)=00g_1(1)0$, but 00 is not a suffix of $g_1(d)$ for any
$d\in\Sigma_4$.
\item
\begin{enumerate}
\renewcommand{\labelenumiii}{(\Alph{enumiii})}
\item Trivial.
\item $g_2(01)=01g_2(0)0$ and $1g_2(01)11=g_2(302)$,
but 11 is not a prefix of $g_2(d)$ for any $d\in\Sigma_4$.
\item $g_1(31)=1g_1(1)01$ and $00g_1(31)0=g_1(012)$,
but 00 is not a suffix of $g_1(d)$ for any $d\in\Sigma_4$.
\item $g_1(23)=0g_1(1)10$ and $01g_1(23)1=g_1(211)$, but $e'=c=1$.
\item $g_1(21)=0g_1(1)01$ and $01g_1(21)0=g_1(212)$, but
$ece'=212\in\mathcal{A}$.
\item $g_2(30)=1g_2(0)10$, $01g_2(30)0=g_2(201)$, and
$g_2(330)0=1g_2(201)$, but $k=a=3$.
\item $g_1(12)=10g_1(1)0$, $0g_1(12)01=g_1(210)$, and
$0g_1(122)=g_1(210)0$, but $k=b=2$.
\end{enumerate}
\end{enumerate}

\begin{center}
\begin{tabular}{|c|c|c|c|}
\hline
Case & $g_i(ab)=tg_i(c)u$ & $vg_i(ab)w=g_i(ece')$ & $g_i(kab)x=yg_i(ece')$ or \\
& & & $xg_i(abk)=g_i(ece')y$ \\
\hline
a.ii & $g_2(01)=0g_2(3)\mathbf{00}$ & - & - \\
     & $g_2(02)=0g_2(1)\mathbf{11}$ & - & - \\
     & $g_2(31)=1g_2(2)\mathbf{00}$ & - & - \\
     & $g_2(32)=1g_2(0)\mathbf{11}$ & - & - \\
\hline
a.iii & $g_1(01)=\mathbf{00}g_1(3)1$ & - & - \\
      & $g_1(02)=\mathbf{00}g_1(1)0$ & - & - \\
      & $g_1(31)=\mathbf{11}g_1(2)1$ & - & - \\
      & $g_1(32)=\mathbf{11}g_1(0)0$ & - & - \\
\hline
a.iv.B & $g_2(01)=01g_2(0)0$ & $1g_2(01)\mathbf{11}=g_2(302)$ & - \\
       & $g_2(32)=10g_2(3)1$ & $0g_2(32)\mathbf{00}=g_2(031)$ & - \\
\hline
a.iv.C & $g_1(02)=0g_1(2)10$ & $\mathbf{11}g_1(02)1=g_1(321)$ & - \\
       & $g_1(31)=1g_1(1)01$ & $\mathbf{00}g_1(31)0=g_1(012)$ & - \\
\hline
\end{tabular}

\medskip

Table 1:  Forbidden Patterns in the Proof of Lemma~\ref{gMainLem}

\end{center}

\begin{center}
\begin{tabular}{|c|c|c|c|}
\hline
Case & $g_i(ab)=tg_i(c)u$ & $vg_i(ab)w=g_i(ece')$ & $g_i(kab)x=yg_i(ece')$ or \\
& & & $xg_i(abk)=g_i(ece')y$ \\
\hline
a.iv.D & $g_1(10)=1g_1(2)01$ & $00g_1(02)0=g_1(0\mathbf{22})$ & - \\
       & $g_1(10)=1g_1(2)01$ & $10g_1(02)0=g_1(1\mathbf{22})$ & - \\
       & $g_1(02)=0g_1(2)10$ & $01g_1(02)1=g_1(\mathbf{22}1)$ & - \\
       & $g_1(23)=0g_1(1)10$ & $01g_1(23)1=g_1(2\mathbf{11})$ & - \\
       & $g_1(23)=0g_1(1)10$ & $11g_1(23)1=g_1(3\mathbf{11})$ & - \\
       & $g_1(31)=1g_1(1)01$ & $10g_1(31)0=g_1(\mathbf{11}2)$ & - \\
       & $g_2(01)=1g_2(2)01$ & $1g_2(01)10=g_2(3\mathbf{00})$ & - \\
       & $g_2(13)=1g_2(2)01$ & $0g_2(13)00=g_2(\mathbf{00}1)$ & - \\
       & $g_2(13)=0g_2(2)10$ & $0g_2(13)01=g_2(\mathbf{00}3)$ & - \\
       & $g_2(20)=0g_2(1)10$ & $1g_2(20)10=g_2(\mathbf{33}0)$ & - \\
       & $g_2(20)=0g_2(1)10$ & $1g_2(20)11=g_2(\mathbf{33}2)$ & - \\
       & $g_2(32)=10g_2(3)1$ & $0g_2(32)01=g_2(0\mathbf{33})$ & - \\
\hline
a.iv.E & $g_1(12)=10g_1(1)0$ & $0g_1(12)10=g_1(\mathbf{212})$ & - \\
       & $g_1(12)=10g_1(1)0$ & $1g_1(12)10=g_1(\mathbf{312})$ & - \\
       & $g_1(12)=1g_1(2)10$ & $00g_1(12)1=g_1(\mathbf{021})$ & - \\
       & $g_1(12)=1g_1(2)10$ & $10g_1(12)1=g_1(\mathbf{121})$ & - \\
       & $g_1(21)=0g_1(1)01$ & $01g_1(21)0=g_1(\mathbf{212})$ & - \\
       & $g_1(21)=0g_1(1)01$ & $11g_1(21)0=g_1(\mathbf{312})$ & - \\
       & $g_1(21)=01g_1(2)1$ & $0g_1(21)01=g_1(\mathbf{021})$ & - \\
       & $g_1(21)=01g_1(2)1$ & $1g_1(21)01=g_1(\mathbf{121})$ & - \\
       & $g_2(03)=01g_2(0)1$ & $1g_2(03)00=g_2(\mathbf{301})$ & - \\
       & $g_2(03)=01g_2(0)1$ & $1g_2(03)01=g_2(\mathbf{303})$ & - \\
       & $g_2(03)=0g_2(3)01$ & $01g_2(03)0=g_2(\mathbf{030})$ & - \\
       & $g_2(03)=0g_2(3)01$ & $01g_2(03)1=g_2(\mathbf{032})$ & - \\
       & $g_2(30)=1g_2(0)10$ & $10g_2(30)0=g_2(\mathbf{301})$ & - \\
       & $g_2(30)=1g_2(0)10$ & $10g_2(30)1=g_2(\mathbf{303})$ & - \\
       & $g_2(30)=10g_2(3)0$ & $0g_2(30)10=g_2(\mathbf{030})$ & - \\
       & $g_2(30)=10g_2(3)0$ & $0g_2(30)11=g_2(\mathbf{032})$ & - \\
\hline
a.iv.F
& $g_2(03)=0g_2(3)01$ & $10g_2(03)0=g_2(130)$ &
  $g_2(\mathbf{00}3)0=0g_2(130)$ \\
& $g_2(03)=0g_2(3)01$ & $10g_2(03)1=g_2(132)$ &
  $g_2(\mathbf{00}3)1=0g_2(132)$ \\
& $g_2(30)=1g_2(0)10$ & $01g_2(30)0=g_2(201)$ &
  $g_2(\mathbf{33}0)0=1g_2(201)$ \\
& $g_2(30)=1g_2(0)10$ & $01g_2(30)1=g_2(203)$ &
  $g_2(\mathbf{33}0)1=1g_2(203)$ \\
\hline
a.iv.G
& $g_1(12)=10g_1(1)0$ & $0g_1(12)01=g_1(210)$ &
  $0g_1(1\mathbf{22})=g_1(210)0$ \\
& $g_1(12)=10g_1(1)0$ & $1g_1(12)01=g_1(310)$ &
  $1g_1(1\mathbf{22})=g_1(310)0$ \\
& $g_1(21)=01g_1(2)1$ & $0g_1(21)10=g_1(023)$ &
  $0g_1(2\mathbf{11})=g_1(023)1$ \\
& $g_1(21)=01g_1(2)1$ & $1g_1(21)10=g_1(123)$ &
  $1g_1(2\mathbf{11})=g_1(123)1$ \\
\hline
\end{tabular}

\medskip

Table 1 (continued): Forbidden Patterns in the Proof of Lemma~\ref{gMainLem}
\end{center}

\item The only $a,b,c$ that satisfy $g_1(a)=st$, $g_1(b)=uv$, and
$g_1(c)=sv$ such that $a\neq c$ and $b\neq c$ are
$(a,b,c)\in\{(0,3,2),(1,2,3),(2,1,0),(3,0,1)\}$.
But by Lemma~\ref{hAvoid}, the infinite word $h^\omega(0)$ contains no subwords of the form
$3\alpha2\alpha0$, $2\alpha3\alpha1$, $1\alpha0\alpha2$,
or $0\alpha1\alpha3$.  This contradicts the assumption that
$b\alpha c\alpha a$ is a subword of $h^\omega(0)$ for some
$\alpha\in\Sigma_4^*$.  The same result holds true for $g_2$.
\end{enumerate}
\endpf
\end{proof}

\begin{lemma}
\label{gSQF}
Neither $g_1(h^\omega(0))$ nor $g_2(h^\omega(0))$ contain squares $yy$
with $|y|\geq 4$.
\end{lemma}

\begin{proof}
As in the case of Lemma~\ref{ming}, this follows from Lemma~\ref{gMainLem}.
\endpf
\end{proof}

\bigskip

   We can now complete the proof of
Theorem~\ref{pu-thm}.  Let ${\bf x} := g_2 (h^\omega(0)) =  
010100011 101 010 011 \cdots$ and
${\bf y} := g_1 (h^\omega(0)) = 
001101010 110 001 010 \cdots$.
Then by Lemma~\ref{Shuf} we have
${\bf x} \sha {\bf y} = f^\omega(0)$.  
But $f^\omega(0)=f^\omega(001)$ and so
$f^\omega(0)$ begins with $f^n(0)f^n(0)$ for all $n\geq0$.
Hence $f^\omega(0)$ begins with an arbitrarily large square.

On the other hand, by
Lemma~\ref{gSQF}, we have that $\bf x$ and $\bf y$ avoid all squares
$ww$ with $|w| \geq 4$.
\endpf
\end{proof}

\section{Acknowledgments}

    We thank Jean-Paul Allouche for helpful discussions.

\end{document}